\documentclass[a4paper,12pt]{amsart}
\usepackage{amsfonts}
\usepackage{amssymb}
\usepackage{ifthen}
\usepackage{graphicx}
\nonstopmode \numberwithin{equation}{section}
\usepackage[usenames]{color}
\newtheorem{thm}{Theorem}
\newtheorem{lem}{Lemma}
\newtheorem{cor}{Corollary}[section]

\newtheorem{cl}{Claim}
\newtheorem{ca}{Case}
\newtheorem{sca}{Subcase}
\newtheorem{scl}{Subclaim}
\newtheorem{conj}[equation]{Conjecture}

\theoremstyle{definition}
\newtheorem{defn}{Definition}

\newtheorem{op}[equation]{Open Problem}
\newtheorem{ques}[equation]{Question}
\newtheorem{rem}{Remark}[section]
\newtheorem{exam}[equation]{Example}

\newcounter {own}
\def\theown {\thesection       .\arabic{own}}

\newenvironment{pf}[1][]{%
 \vskip 3mm
 \noindent
 \ifthenelse{\equal{#1}{}}%
  {{\slshape Proof. }}%
  {{\slshape #1.} }%
 }%
{\qed\bigskip}

\newcounter{alphabet}
\newcounter{tmp}
\newenvironment{Thm}[1][]{\refstepcounter{alphabet}%
\bigskip%
\noindent%
{\bf Theorem \Alph{alphabet}}%
\ifthenelse{\equal{#1}{}}{}{ (#1)}%
{\bf .} \itshape}{\vskip 8pt}

\makeatletter
\newcommand{\Ref}[1]{\@ifundefined{r@#1}{}{\setcounter{tmp}{\ref{#1}}\Alph{tmp}}}
\makeatother

\newcommand{\IC}{{\mathbb C}}
\newcommand{\ID}{{\mathbb D}}
\newcommand{\IB}{{\mathbb B}}

\def\be{\begin{equation}}
\def\ee{\end{equation}}

\newcommand{\bee}{\begin{enumerate}}
\newcommand{\eee}{\end{enumerate}}

\newcommand{\blem}{\begin{lem}}
\newcommand{\elem}{\end{lem}}
\newcommand{\bthm}{\begin{thm}}
\newcommand{\ethm}{\end{thm}}
\newcommand{\bcor}{\begin{cor}}
\newcommand{\ecor}{\end{cor}}
\newcommand{\beg}{\begin{exam}}
\newcommand{\eeg}{\end{exam}}
\newcommand{\begs}{\begin{examples}}
\newcommand{\eegs}{\end{examples}}
\newcommand{\bdefe}{\begin{defn}}
\newcommand{\edefe}{\end{defn}}
\newcommand{\bprob}{\begin{prob}}
\newcommand{\eprob}{\end{prob}}
\newcommand{\bques}{\begin{ques}}
\newcommand{\eques}{\end{ques}}
\newcommand{\bei}{\begin{itemize}}
\newcommand{\eei}{\end{itemize}}
\newcommand{\bcon}{\begin{conj}}
\newcommand{\econ}{\end{conj}}
\newcommand{\bop}{\begin{op}}
\newcommand{\eop}{\end{op}}

\newcommand{\bca}{\begin{ca}}
\newcommand{\eca}{\end{ca}}
\newcommand{\bsca}{\begin{sca}}
\newcommand{\esca}{\end{sca}}

\newcommand{\bcl}{\begin{cl}}
\newcommand{\ecl}{\end{cl}}

\newcommand{\bscl}{\begin{scl}}
\newcommand{\escl}{\end{scl}}

\newcommand{\bcons}{\begin{conjs}}
\newcommand{\econs}{\end{conjs}}
\newcommand{\bprop}{\begin{propo}}
\newcommand{\eprop}{\end{propo}}
\newcommand{\br}{\begin{rem}}
\newcommand{\er}{\end{rem}}
\newcommand{\brs}{\begin{rems}}
\newcommand{\ers}{\end{rems}}
\newcommand{\bo}{\begin{obser}}
\newcommand{\eo}{\end{obser}}
\newcommand{\bos}{\begin{obsers}}
\newcommand{\eos}{\end{obsers}}
\newcommand{\bpf}{\begin{pf}}
\newcommand{\epf}{\end{pf}}
\newcommand{\ba}{\begin{array}}
\newcommand{\ea}{\end{array}}
\newcommand{\beq}{\begin{eqnarray}}
\newcommand{\beqq}{\begin{eqnarray*}}
\newcommand{\eeq}{\end{eqnarray}}
\newcommand{\eeqq}{\end{eqnarray*}}

\newcommand{\ds}{\displaystyle}

\newcounter{minutes}
\divide\time by 60
\newcounter{hours}
\multiply\time by 60 \addtocounter{minutes}{-\time}

\begin{document}
\bibliographystyle{amsplain}
\title [] {Radial growth, Lipschitz  and Dirichlet spaces on solutions to the Yukawa equation}

\def\thefootnote{}
\footnotetext{ \texttt{\tiny File:~\jobname .tex,
          printed: \number\day-\number\month-\number\year,
          \thehours.\ifnum\theminutes<10{0}\fi\theminutes}
} \makeatletter\def\thefootnote{\@arabic\c@footnote}\makeatother


\author{SH. Chen}
\address{Sh. Chen, Department of Mathematics,
Hunan Normal University, Changsha, Hunan 410081, People's Republic
of China.} \email{shlchen1982@yahoo.com.cn}

\author{ A. Rasila }
\address{A. Rasila, Department of Mathematics and Systems Analysis, Aalto University, P. O. Box 11100, FI-00076 Aalto,
 Finland.} \email{antti.rasila@iki.fi}

\author{X. Wang${}^{~\mathbf{*}}$}
\address{X. Wang, Department of Mathematics,
Hunan Normal University, Changsha, Hunan 410081, People's Republic
of China.} \email{xtwang@hunnu.edu.cn}

\keywords{Yukawa  PDE, Green's theorem, Lipschitz space.\\
${}^{\mathbf{*}}$ Corresponding author}

\begin{abstract}
In this paper, we investigate some properties to  solutions $f$ to
the Yukawa PDE: $\Delta f=\lambda f$ in the unit ball $\IB^n$ of
$\IC^n$, where $\lambda$ is a nonnegative constant. First, we prove
that the answer to an open problem of Girela and Pel\'{a}ez,
concerning such solutions, is positive. Then we study  relationships
on such solutions between the bounded mean oscillation and
Lipschitz-type spaces. At last, we discuss Dirichlet-type energy
integrals on such solutions in the unit ball of $\mathbb{C}^{n}$ and
give an application.
\end{abstract}

\thanks{The research was partly supported by
NSF of China (No. 11071063)} 

\maketitle \pagestyle{myheadings} \markboth{Sh. Chen,  A. Rasila and
X. Wang}{Radial growth, Lipschitz  and Dirichlet spaces on solutions
to the Yukawa equation}

\section{Introduction and Main results}\label{csw-sec1}
Let $\mathbb{C}$ denote the complex plane. We write
$\mathbb{C}^{n}=\{z=(z_{1},\ldots,z_{n}):\,
z_{1},\ldots,z_{n}\in\mathbb{C}\}$, $$\IB^n(a,
r)=\left\{z\in\mathbb{C}^{n}:\,
|z-a|=\Big(\sum^{n}_{k=1}|z_k-a_k|^2\Big)^{1/2}<1\right\} $$ and
$\IB^n=\IB^n(0, 1)$, the unit ball in $\mathbb{C}^{n}$. In
particular, the unit disk of $\mathbb{C}$, i.e. $\IB^2$, is denoted
by $\mathbb{D}$. We use $d(z)$ to denote  the Euclidean distance
from $z$ to the boundary  of $\mathbb{B}^{n}$. Let $\lambda$ be a
nonnegative constant and $f=u+iv$ be a complex-valued function of
$\mathbb{B}^{n}$ into $\mathbb{C}$, where $u$ and $v$ are
real-valued and twice continuously differentiable functions of
$\mathbb{B}^{n}$ into $\mathbb{R}$. The following elliptic partial
differential equation, or briefly PDE in the following,
\be\label{eq-g2} \Delta f(z)=\lambda f(z)\ee in $\IB^n$ is called
the {\it Yukawa PDE}, where $\Delta$ represents the usual complex
Laplacian operator
$$\Delta:= \sum_{k=1}^{n}\left(\frac{\partial^{2}}{\partial
x^{2}_{k}}+\frac{\partial^{2}}{\partial
y^{2}_{k}}\right)=4\sum_{k=1}^{n}\frac{\partial^{2}}{\partial
z_{k}\partial\overline{z}_{k}}
$$
for each $k\in\{1,\ldots,n\}$. Here $z_{k}=x_{k}+iy_{k}$.

Equation (\ref{eq-g2}) arose out of an attempt by the Japanese
physicist Hideki Yukawa to describe the nuclear potential of a point
charge as $e^{-\sqrt{\lambda} r}/r$ (cf. \cite{D-,SW}). It is well
known that each solution $f$ to (\ref{eq-g2}) belongs to
$C^{\infty}(\mathbb{B}^{n})$, i.e., they are infinitely
differentiable in $\mathbb{B}^{n}$. We refer to \cite{A,B,Ya} for
basic results on the theory of elliptic PDEs. Moreover, if
$\lambda\leq0$ in (\ref{eq-g2}), then  (\ref{eq-g2}) is called {\it
Helmholtz equation} (see \cite{E}). Especially, if $\lambda=0$ in
(\ref{eq-g2}), then $f$ is a complex-valued {\it harmonic} mapping
(cf. \cite{ABR}).  Moreover, if $\lambda=0$ in (\ref{eq-g2}) with
$n=1$, then $f$ is a complex-valued planar harmonic mapping. It is
known that every planar harmonic mapping $f$ defined in  $\ID$
admits a decomposition $f=h+\overline{g}$, where $h$ and $g$ are
analytic in $\ID$. We refer to \cite{Du} for basic results
concerning planar harmonic mappings.

For a complex-valued and differentiable function  $f$ of
$\mathbb{B}^{n}$ into $\mathbb{C}$, we introduce the following
notations (cf. \cite{CPW1-1,CPW1-2}):
$$f_{z}=(f_{z_{1}},\ldots,f_{z_{n}}),~
f_{\overline{z}}=(f_{\overline{z}_{1}},\ldots,f_{\overline{z}_{n}})~\mbox{and}~\widetilde{\nabla
f}=(f_{z},f_{\overline{z}}). $$ Let
 $|\widetilde{\nabla f}|$ be the {\it Hilbert-Schmidt norm} given by $$|\widetilde{\nabla
f}|=(|f_{z}|^{2}+|f_{\overline{z}}|^{2})^{1/2}.$$

Let $f=u+iv$ be a continuously differentiable mapping from
$\mathbb{B}^{n}$ into $\mathbb{C}$, where $u$ and $v$ are
real-valued functions. Then for
$z=(z_{1},\cdots,z_{n})=(x_{1}+iy_{1},\cdots,x_{n}+iy_{n})\in\mathbb{B}^{n}$,
\be\label{1x} |f_{z}(z)|+|f_{\overline{z}}(z)|\leq |\nabla
u(z)|+|\nabla v(z)|, \ee where $\nabla u=\Big( \frac{\partial
u}{\partial x_{1}}, \frac{\partial u}{\partial y_{1}},\cdots,
\frac{\partial u}{\partial x_{n}}, \frac{\partial u}{\partial
y_{n}}\Big)$ and $\nabla v=\Big( \frac{\partial v}{\partial x_{1}},
\frac{\partial v}{\partial y_{1}},\cdots, \frac{\partial v}{\partial
x_{n}}, \frac{\partial v}{\partial y_{n}} \Big)$. But the converse
of (\ref{1x}) is not always true (see \cite{CPW3}).

For $p\in(0,\infty]$, the Hardy space $\mathcal{H}^{p}$
 consists of those functions $f:\ \mathbb{B}^{n}\rightarrow\mathbb{C}$
 such that $f$ is measurable, $M_{p}(r,f)$ exists for all $r\in(0,1)$ and  $ \|f\|_{p}<\infty$, where
$$
\|f\|_{p}=
\begin{cases}
\displaystyle\sup_{0<r<1}M_{p}(r,f),
& \mbox{ if } p\in(0,\infty),\\
\displaystyle\sup_{z\in\mathbb{B}^{n}}|f(z)|, &\mbox{ if }\,
p=\infty,
\end{cases}~ M_{p}(r,f)=\left(\int_{\partial\mathbb{B}^{n}}|f(r\zeta)|^{p}\,d\sigma(\zeta)\right)^{1/p}
$$
and $d\sigma$ denotes the normalized Lebesgue surface measure in
$\partial\mathbb{B}^{n}$.

A continuous increasing function $\omega:\, [0,+\infty)\rightarrow
[0,+\infty)$ with $\omega(0)=0$ is called a {\it majorant} if
$\omega(t)/t$ is non-increasing for $t>0$. Given a subset $\Omega$
of $\mathbb{C}$, a function $f:\, \Omega\rightarrow \mathbb{C}$ is
said to belong to the {\it Lipschitz space
$\Lambda_{\omega}(\Omega)$} if there is a positive constant $C$ such
that \be\label{eq1} |f(z)-f(w)|\leq C\omega(|z-w|) ~\mbox{ for all
$z,\ w\in\Omega.$} \ee For $\delta_{0}>0$, let \be\label{eq2}
\int_{0}^{\delta}\frac{\omega(t)}{t}\,dt\leq C\cdot\omega(\delta),\
0<\delta<\delta_{0} \ee and \be\label{eq3}
\delta\int_{\delta}^{+\infty}\frac{\omega(t)}{t^{2}}\,dt\leq
C\cdot\omega(\delta),\ 0<\delta<\delta_{0}, \ee where $\omega$ is a
majorant and $C$ is a positive constant. A majorant $\omega$ is said
to be {\it regular} if it satisfies the conditions (\ref{eq2}) and
(\ref{eq3}) (see \cite{D,P}).

In \cite{GP}, Girela and Pel\'{a}ez obtained the following result.

\begin{Thm} $($\cite[Theorem 1(a) ]{GP}$)$
\label{ThmB} Let $p\in(2,\infty)$. For $r\in(0,1)$, if  $f$ is an
analytic function in $\mathbb{D}$ such that
$$M_{p}(r,f')=O \left (\Big(\frac{1}{1-r}\Big)\right  )  ~\mbox{ as
$r\rightarrow1$},
$$
then for all $\ds \beta>1/2$,

\be\label{eq1.1a} M_{p}(r,f)=O \left
(\Big(\log\frac{1}{1-r}\Big)^{\beta} \right ) ~\mbox{ as
$r\rightarrow 1$}.  \ee
\end{Thm}

In \cite[P$_{464}$, Equation (26)]{GP}, Girela and Pel\'{a}ez asked
whether $\beta$ in (\ref{eq1.1a}) can be substituted by $1/2$. This
problem was affirmatively settled by Girela, Pavlovic and Pel\'{a}ez
in \cite{GPP}. In \cite{CPW2}, the authors proved further that the
answer to this problem is affirmative for the setting of
complex-valued harmonic mappings in $\ID$. The first aim of this
paper is to show that the answer to this problem is also affirmative
for mappings $f$ satisfying (\ref{eq-g2}) and $p\in[2,\infty)$. Our
result is given as follows.

\begin{thm}\label{thm3.2}
Let  $p\in[2,\infty)$,   $\lambda \in[0,4n/p)$ and $\omega$ be a
majorant. For $r\in(0,1)$, if $f$ is a solution to \eqref{eq-g2}
such that
$$ M_{p}(r,\widetilde{\nabla f})\leq C\omega\Big(\frac{1}{1-r}\Big), $$ then

$$M_{p}(r,f)\leq\Big(\frac{4n}{4n-p\lambda}\Big)^{1/2}\Big(|f(0)|^{2}+2p(p-1)C^{2}\omega(1)T(r)\Big)^{1/2},$$ where
\[
T(r)=\int_{0}^{1}\omega\left(\frac{1}{1-\rho r}\right)d\rho
\]
and $C$ is a positive constant.
\end{thm}

By taking $\omega(t)=t$ in Theorem \ref{thm3.2}, we obtain the
following result.
\begin{cor}\label{thm1-1}
Let $p\in[2,\infty)$ and  $\lambda \in[0,4n/p)$. For $r\in(0,1)$, if
$f$ is a solution to \eqref{eq-g2} such that
$$ M_{p}(r,\widetilde{\nabla f})=O\left (\Big(\frac{1}{1-r}\Big)\right ) ~\mbox{
as $r\rightarrow 1$}, $$ then

$$M_{p}(r,f)=O \left (\Big(\log\frac{1}{1-r}\Big)^{1/2} \right )~\mbox{
as}~ r\rightarrow1.$$
\end{cor}

\begin{rem} Obviously, all analytic functions
and complex-valued harmonic mappings defined in $\mathbb{B}^{n}$ are
solutions to \eqref{eq-g2} with $\lambda=0$, and there also are
solutions which are neither analytic nor harmonic. For example, we
can take $f(z)=e^{[\sum_{k=1}^{n}(z_{k}+\overline{z}_{k}/2)]},$
where $z\in\mathbb{B}^{n}$.
Hence Theorem \ref{thm3.2} and Corollary \ref{thm1-1} are
generalizations of \cite[Theorem 1.1]{GPP}, \cite[Theorem
1(a)]{CPW2} and \cite[Corollary 6]{ST}. But it is not clear for us
that what the best upper bound of $\lambda$ in Theorem \ref{thm3.2}
is.
\end{rem}

In \cite{Kr}, the author discussed the Lipschitz spaces on smooth
functions. Dyakonov \cite{D} discussed the relationship between the
Lipschitz space and the bounded mean oscillation on analytic
functions in $\mathbb{D}$, and obtained the following result.

\begin{Thm}{\rm \cite[Theorem 1]{D}}\label{ThmA}
Suppose that $f$ is a analytic fucntion in $\mathbb{D}$ which is
continuous up to the boundary of $\mathbb{D}$. If $\omega$ and
$\omega^{2}$ are regular majorants,  then
$$f\in L_{\omega}(\mathbb{D})\Longleftrightarrow \left(\mathcal{P}_{|f|^{2}}(z)-|f(z)|^{2}\right)^{1/2}\leq C\omega(d(z)),$$
where
$\mathcal{P}_{|f|^{2}}(z)=\frac{1}{2\pi}\int_{0}^{2\pi}\frac{1-|z|^{2}}{|z-e^{i\theta}|^{2}}|f(e^{i\theta})|^{2}d\theta$
and $C$ is a positive constant.
\end{Thm}

For the solutions to \eqref{eq-g2}, 
we also get the following theorem, which is similar to Theorem
\Ref{ThmA}.

\begin{thm}\label{thm1-y}
Let $\omega$ be a majorant and $f$ be a solution to \eqref{eq-g2}.
If $f$ satisfies
$$|\widetilde{\nabla f}(z)|\leq C\omega\Big(\frac{1}{d(z)}\Big)$$ in
$\mathbb{B}^{n}$, then for all $r\in(0,d(z)]$,
$$\frac{1}{|\mathbb{B}^{n}(z,r)|}\int_{\mathbb{B}^{n}(z,r)}|f(\zeta)-f(z)|dV(\zeta)\leq Cr\omega\Big(\frac{1}{r}\Big),$$ where
 $C$ is a positive constant and
$dV$ denotes the Lebesgue volume measure in $\mathbb{B}^{n}$.
\end{thm}

In particular, if $f$  is a solution to \eqref{eq-g2} with
$\lambda=0$, then we have

\begin{thm}\label{thm1y}
Let $\omega$ be a majorant and $f$ be a solution to \eqref{eq-g2}
with $\lambda=0$. Then $f$ satisfies
$$|\widetilde{\nabla f}(z)|\leq C\omega\Big(\frac{1}{d(z)}\Big)$$ in $\mathbb{B}^{n}$
if and only if for all $r\in(0,d(z)]$,
$$\frac{1}{|\mathbb{B}^{n}(z,r)|}\int_{\mathbb{B}^{n}(z,r)}|f(\zeta)-f(z)|dV(\zeta)\leq Cr\omega\Big(\frac{1}{r}\Big),$$ where
 $C$ is a positive constant.
\end{thm}

\begin{defn}\label{defn2}
Let $f$ be a continuous function in $\mathbb{B}^{n}$. We say $f\in
BMO$ if
$$\|f\|_{BMO}=\sup_{\mathbb{B}^{n}(z,r)\subseteq\mathbb{B}^{n}}\frac{1}{|\mathbb{B}^{n}(z,r)|}\int_{\mathbb{B}^{n}(z,r)}\left|f(\zeta)-
\frac{1}{|\mathbb{B}^{n}(z,r)|}\int_{\mathbb{B}^{n}(z,r)}f(\xi)dV(\xi)\right|dV(\zeta)
$$ is bounded, where   $r\in(0,d(z)]$.
\end{defn}

In particular, by taking $\omega(t)=t$ in Theorem \ref{thm1y}, we
get the following result.
\begin{cor}\label{cor-1}
Let $f$ be a solution to \eqref{eq-g2} with $\lambda=0$. Then $f\in
BMO$ if and only if $|\widetilde{\nabla f}(z)|\leq M\frac{1}{d(z)}$
holds in $\mathbb{B}^{n}$.
\end{cor}

For $\nu,~\gamma,~t\in\mathbb{R}$,
 $$ D_{f}(\nu,\gamma,t)=\int_{\mathbb{B}^{n}}(1-|z|)^{\nu}|f(z)|^{\gamma}|\widetilde{\nabla
 f}(z)|^{t}dV_{N}(z)$$
is called {\it  Dirichlet-type energy integral } of $f$ defined in
$\mathbb{B}^{n}$, where $dV_{N}$ denotes the normalized  Lebesgue
volume measure in $\mathbb{B}^{n}$ (cf. \cite{E, SH,ST,W}).



\begin{thm}\label{thm4}
Let $f$ be a solution to \eqref{eq-g2}. Then there exist positive
constants $C_{1}$ and $C_{2}$ such that
$$\int_{\mathbb{D}}(1-|z|)^{1+\frac{2}{\beta}(n-1)}\Delta(|f(z)|^{\frac{2}{\beta}})dV_{N}(z)\leq C_{1}D_{f}(\beta-1,1,1)+C_{2},$$ where
$\beta\in(0,1]$.
\end{thm}

As an application of Theorem \ref{thm4}, we get the following
result.

\begin{cor}\label{cor-1}
Let $f$ be a solution to \eqref{eq-g2}. If $n=1$ and
$D_{f}(\beta-1,1,1)<\infty$, then
$f\in\mathcal{H}^{\frac{2}{\beta}}$, where $\beta\in(0,1]$.
\end{cor}

\section{Integral means and Lipschitz spaces}\label{csw-sec2}

We start this section by recalling the following result  $($cf.
\cite{Pav,ST,Z}$)$.


\begin{Thm}{\rm \bf (Green's Theorem)}\label{Green-thm}
Let $g$ be a function of class $C^{2}(\mathbb{B}^{n})$. If $n\geq2$,
then for $r\in(0,1)$,
$$\int_{\partial \mathbb{B}^{n}}g(r\zeta)\,d\sigma(\zeta)=
g(0)+\int_{ \mathbb{B}^{n}(0,r)}\Delta g(z)G_{2n}(z,r)\,dV_{N}(z),
$$
where $G_{2n}(z,r)=(|z|^{2(1-n)}-r^{2(1-n)})/[4n(n-1)]$. Moreover,
if $n=1$, then for $r\in(0,1)$,
\begin{eqnarray*}
\frac{1}{2\pi}\int_{0}^{2\pi}g(re^{i\theta})\,d\theta
&=&g(0)+ \frac{1}{2}\int_{\mathbb{D}_{r}}\Delta
g(z)\log\frac{r}{|z|}\,dA(z),
\end{eqnarray*}  where $dA$ denotes the normalized area measure in
$\mathbb{D}$.
\end{Thm}



Recall that a real-valued and   continuous function $u$ defined in
 $\mathbb{B}^{n}$ is {\it subharmonic} if for all $z_{0}\in
\mathbb{B}^{n}$, there is $\varepsilon\in(0,1-|z_{0}|)$ such that
$$u(z_{0})\leq\int_{\partial\mathbb{B}^{n}}u(z_{0}+r\zeta)d\sigma(\zeta)$$ holds
for all $ r\in[0,\varepsilon)$. Moreover, if $u\in
C^{2}(\mathbb{B}^{n})$, then $u$ is subharmonic if and only if
$\Delta u\geq0$ in $\mathbb{B}^{n}$ (cf. \cite{Do}).


\begin{lem}\label{lem1}
Suppose that $f$ is a solution to \eqref{eq-g2}. Then

\noindent  $(I)$ for $p\in[2,\infty)$, $M_{p}^{p}(r,f)$ is
increasing  in $(0,1)$ and $| f|^{p}$ is subharmonic in
$\mathbb{B}^{n}$;

\noindent $(II)$ $M_{2}^{2}(r,\widetilde{\nabla f})$ is increasing
 in $(0,1)$ and $| \widetilde{\nabla f}|^{2}$ is subharmonic
in $\mathbb{B}^{n}$;

Moreover, if $f$ is a solution of \eqref{eq-g2} with $\lambda=0$,
then $|f|^{p}$ is subharmonic in $\mathbb{B}^{n}$ for
$p\in[1,\infty)$.
\end{lem}
\bpf We first prove $(I)$. For this, we consider the case where
$p\in[2,4)$ and the case where $p\in[4,\infty)$, separately. \bca
Suppose first $p\in[2,4)$.\eca  Let
$F_{m}^{p}=(|f|^{2}+\frac{1}{m})^{p/2}$. By elementary calculations,
we have
\begin{eqnarray*}
\Delta(F_{m}^{p})&=&4\sum_{k=1}^{n}\frac{\partial^{2}}{\partial
z_{k}\partial\overline{z}_{k}}(F_{m}^{p})\\
&=&p(p-2)(|f|^{2}+\frac{1}{m})^{\frac{p}{2}-2}\sum_{k=1}^{n}\left|f\overline{f}_{z_{k}}+f_{\overline{z}_{k}}\overline{f}\right|^{2}\\
&&+2p(|f|^{2} +\frac{1}{m})^{\frac{p}{2}-1}|\widetilde{\nabla
f}|^{2}+p\lambda|f|^{2}(|f|^{2}+\frac{1}{m})^{\frac{p}{2}-1}.
\end{eqnarray*}
Let $T_{m}=\Delta(F_{m}^{p})$. Obviously, for   $r\in(0,1)$, $T_{m}$
is integrable in $\mathbb{B}^{n}(0,r)$ and $T_{m}\leq F,$ where
$$F=p(p-2)|f|^{p-2}\sum_{k=1}^{n}(|f_{z_{k}}|+|f_{\overline{z}_{k}}|)^{2}+2p(1+|f|^{2})^{\frac{p}{2}-1}|\widetilde{\nabla
f}|^{2}+p\lambda|f|^{2}(|f|^{2}+1)^{\frac{p}{2}-1}$$ and $F$ is
integrable in $\mathbb{B}^{n}(0,r)$. By Theorem \Ref{Green-thm} and
Lebesgue's dominated convergence Theorem, we have
\begin{eqnarray*}
\lim_{m\rightarrow\infty}r^{2n-1}\frac{d}{dr}M_{p}^{p}(r,F_{m})&=&\frac{1}{2n}
\lim_{m\rightarrow\infty}\int_{\mathbb{B}^{n}(0,r)}T_{m}dV_{N}(z)\\
&=&\frac{1}{2n}\int_{\mathbb{B}^{n}(0,r)}\lim_{n\rightarrow\infty}T_{n}dV_{N}(z)\\
&=&\frac{1}{2n}\int_{\mathbb{B}^{n}(0,r)}\big[p(p-2)|f|^{p-4}\sum_{k=1}^{n}|f\overline{f}_{z_{k}}+\overline{f}f_{\overline{z}_{k}}|^{2}\\
&&+2p|f|^{p-2}|\widetilde{\nabla
f}|^{2}+p\lambda|f|^{p}\big]dV_{N}(z)\\
 &=&r^{2n-1}\frac{d}{dr}M_{p}^{p}(r,f)\\
 &\geq&0,
\end{eqnarray*}
which implies that $M_{p}^{p}(r,f)$ is increasing   in $(0,1)$  for
$p\in[2,4)$.

\bca\label{case-2} Suppose then $p\in[4,\infty)$.\eca  By
computations, we get
$$\Delta(|f|^{p})=p(p-2)|f|^{p-4}\sum_{k=1}^{n}|f\overline{f}_{z_{k}}+\overline{f}f_{\overline{z}_{k}}|^{2}
+2p|f|^{p-2}|\widetilde{\nabla f}|^{2}+p\lambda|f|^{p}\geq0,$$ which
gives that $M_{p}^{p}(r,f)$ is increasing  in $(0,1)$. 

By Cases 1 and 2, we see that for $p\in[2,\infty)$, $M_{p}^{p}(r,f)$
is increasing  in $(0,1)$. This shows that for every point $z_{0}\in
\mathbb{B}^{n}$,
$$|f(z_{0})|^{p}\leq\int_{\partial\mathbb{B}^{n}}|f(z_{0}+ r
\zeta)|^{p}d\sigma(\zeta)$$  for all $r\in[0,1-|z_{0}|)$. Hence
$|f|^{p}$ is subharmonic in $\mathbb{B}^{n}$. The proof of $(I)$ is
complete.


Now we come to prove $(II)$.
\begin{eqnarray*}
\Delta(\widetilde{|\nabla
f}|^{2})&=&\Delta\Big[\sum_{k=1}^{n}(f_{z_{k}}\overline{f}_{z_{k}}+f_{\overline{z}_{k}}\overline{f}_{\overline{z}_{k}})\Big]\\
&=&4\lambda|\widetilde{\nabla f}|^{2}+4\sum_{k=1}^{n}\sum_{j=1}^{n}
\left(|f_{z_{k}z_{j}}|^{2}+|f_{\bar{z}_{k} \bar{z}_{j}}|^{2}+|f_{\overline{z}_{k}z_{j}}|^{2}+|f_{z_{k}\overline{z}_{j}}|^{2}\right)\\
&\geq&0,
\end{eqnarray*}
which implies that $M_{2}^{2}(r,\widetilde{\nabla f})$ is increasing
in $(0,1)$ and $|\widetilde{\nabla f}|^{2}$ is subharmonic in
$\mathbb{B}^{n}$.

In particular, if $f$ is a solution to (\ref{eq-g2}) with
$\lambda=0$, then $f$ is a  harmonic mapping. This implies that
$|f|^{p}$ is subharmonic in $\mathbb{B}^{n}$ for $p\in[1,\infty)$
(cf. \cite{ABR}). The proof of this lemma is complete. \epf

By using Theorem \Ref{Green-thm} and the similar argument as in the
proof of Case 1 of Lemma \ref{lem1}, we obtain the following result.

\begin{lem}\label{Lemma-2.0}
Let $p\in[2,\infty)$, $r\in(0,1)$, and suppose that $f$ is a
solution to \eqref{eq-g2}. Then
$$M_{p}^{p}(r,f)=|f(0)|^{p}+ \int_{\mathbb{B}^{n}(0,r)}\Delta
(|f(z)|^{p})G_{2n}(z,r)dV_{N}(z)$$ and
$$r^{2n-1}\frac{d}{dr}M_{p}^{p}(r,f)=\frac{1}{2n}\int_{\mathbb{B}^{n}(0,r)}\Delta\big
(|f(z)|^{p}\big )dV_{N}(z),
 $$ where $G_{2n}$ is the function defined in Theorem {\rm \Ref{Green-thm}}.
\end{lem}

The following result is useful to the proof of Theorem \ref{thm3.2}.

\begin{lem}\label{lem0.1}
Let $p\in[2,\infty)$, $r\in(0,1)$ and  $f$ be a solution to
\eqref{eq-g2}. Then
$$\int_{\mathbb{B}^{n}(0,r)}|f(z)|^{p}G_{2n}(z,r)dV_{N}(z)\leq\frac{r^{2}}{4n}M_{p}^{p}(r,f).$$
\end{lem}
\bpf By Lemma \ref{lem1}, we see that $M_{p}^{p}(\rho,f)$ is
increasing on $\rho\in(0,r]$.  Let
$$I(r)=\int_{\mathbb{B}^{n}(0,r)}|f(z)|^{p}G_{2n}(z,r)dV_{N}(z).$$
Then
\begin{eqnarray*}
I(r)&=&\frac{1}{2(n-1)}\int_{0}^{r}\left[\int_{\partial\mathbb{B}^{n}}
|f(\rho\zeta)|^{p}\Big(\rho-\rho^{2n-1}r^{2(1-n)}\Big)d\sigma(\zeta)\right]d\rho\\
&=&\frac{1}{2(n-1)}\int_{0}^{r}M_{p}^{p}(\rho,f)
\big(\rho-\rho^{2n-1}r^{2(1-n)}\big)d\rho\\
&\leq&\frac{M_{p}^{p}(r,f)}{2(n-1)}\int_{0}^{r}
\big(\rho-\rho^{2n-1}r^{2(1-n)}\big)d\rho\\
&=&\frac{r^{2}}{4n}M_{p}^{p}(r,f).
\end{eqnarray*}
The proof of this lemma is complete. \epf

Now we are ready to prove Theorems \ref{thm3.2} and \ref{thm1-y}.

\subsection*{Proof of  Theorem \ref{thm3.2}} Set
$$\mathcal{A}(r,f)=\int_{\partial\mathbb{B}^{n}}|f(r\zeta)|^{p-2}|\widetilde{\nabla f}(r\zeta)|^{2}\,d\sigma(\zeta).
$$
Then H\"older's inequality yields
\begin{eqnarray*}
\mathcal{A}(r,f)&\leq&\left(\int_{\partial\mathbb{B}^{n}}|\widetilde{\nabla
f}(r\zeta)|^{p}\,d\sigma(\zeta)\right)^{2/p}
\left(\int_{\partial\mathbb{B}^{n}}| f(r\zeta)|^{p}\, d\sigma(\zeta)\right)^{(p-2)/p}\\
&=&M_{p}^{2}(r,\widetilde{\nabla f})\cdot M_{p}^{p-2}(r, f).
\end{eqnarray*}
By using polar coordinates, we see from Lemmas \ref{Lemma-2.0} and
\ref{lem0.1} that
\begin{eqnarray*}
M_{p}^{p}(r, f)&=&|f(0)|^{p}+\int_{\mathbb{B}^{n}(0,r)}\Delta(|f(z)|^{p})G_{2n}(z,r)\,dV_{N}(z)\\
&\leq&|f(0)|^{p}+\int_{\mathbb{B}^{n}(0,r)}\big[2p(p-1)|f(z)|^{p-2}|\widetilde{\nabla
f}(z)|^{2}+\lambda p|f(z)|^{p}\big]
G_{2n}(z,r)\,dV_{N}(z)\\
&=&|f(0)|^{p}+\int_{0}^{r}\int_{\partial\mathbb{B}^{n}}4np(p-1)\rho^{2n-1}|f(\rho\zeta)|^{p-2}|\widetilde{\nabla
f}(\rho\zeta)|^{2}G_{2n}(\rho \zeta,r)d\sigma(\zeta)d\rho\\
&&+p\lambda\int_{\mathbb{B}^{n}(0,r)}|f(z)|^{p}G_{2n}(z,r)\,dV_{N}(z)\\
&=&|f(0)|^{p}+\int_{0}^{r}4np(p-1)\rho^{2n-1}G_{2n}(\rho\zeta,r)\mathcal{A}(\rho,f)d\rho\\
&&+
p\lambda\int_{\mathbb{B}^{n}(0,r)}|f(z)|^{p}G_{2n}(z,r)\,dV_{N}(z)\\
&\leq&|f(0)|^{p}+4p(p-1)\int_{0}^{r}n\rho^{2n-1}G_{2n}(\rho\zeta,r)M_{p}^{2}(\rho,\widetilde{\nabla
f}) M_{p}^{p-2}(\rho, f)\,d\rho\\
&&+\frac{p\lambda r^{2}}{4n}M_{p}^{p}(r,f),
\end{eqnarray*}
which, because $M_p(r,f)$ is increasing on $r$, implies
\begin{eqnarray*}
 \left(1-\frac{p\lambda }{4n}\right)M_{p}^{2}(r, f)&\leq&\left(1-\frac{p\lambda r^{2}}{4n}\right)M_{p}^{2}(r,
 f)\\
 &\leq&|f(0)|^{2}+4p(p-1)\int_{0}^{r}n\rho^{2n-1}G_{2n}(\rho\zeta,r)M_{p}^{2}(\rho,\widetilde{\nabla f})\,d\rho\\
&=&|f(0)|^{2}+2p(p-1)\int_{0}^{1}r^{2} M_{p}^{2}(r\rho,\widetilde{\nabla f})\cdot\frac{\rho(1-\rho^{2n-2})}{2(n-1)}\,d\rho\\
\end{eqnarray*}
\begin{eqnarray*}
&\leq&|f(0)|^{2}+2p(p-1)\int_{0}^{1}M_{p}^{2}(r\rho,\widetilde{\nabla f})(1-\rho)\,d\rho\\
&\leq&|f(0)|^{2}+2p(p-1)C^{2}\int_{0}^{1}\left[\omega\Big(\frac{1}{1-r\rho}\Big)\right]^{2}(1-\rho)\,d\rho\\
&\leq&|f(0)|^{2}+2p(p-1)C^{2}\int_{0}^{1}\left[\omega\Big(\frac{1}{1-r\rho}\Big)\right]^{2}(1-r\rho)\frac{(1-\rho)}{(1-r\rho)}\,d\rho\\
&\leq&|f(0)|^{2}+2p(p-1)C^{2}\omega(1)\int_{0}^{1}\omega\Big(\frac{1}{1-r\rho}\Big)d\rho\\
&=&|f(0)|^{2}+2p(p-1)C^{2}\omega(1)T(r),
\end{eqnarray*}
where $C$ is a positive constant. This observation gives the desired
result:
$$M_{p}(r,f)\leq\Big(\frac{4n}{4n-p\lambda}\Big)^{1/2}\Big(|f(0)|^{2}+2p(p-1)C^{2}\omega(1)T(r)\Big)^{1/2}.$$
The proof of this theorem is complete. \qed

\subsection*{Proof of Theorem \ref{thm1-y}}
For $z,\ w\in\mathbb{B}^{n}$ and $t\in[0,1]$, we have
$$d\big(z+t(w-z)\big)=1-|z+t(w-z)|\geq d(z)-t|w-z|.$$ Suppose that
$d(z)-t|w-z|>0$. Then

\begin{eqnarray*}
|f(z)-f(w)|&=&\Big|\int_{0}^{1}\frac{df}{dt}(wt+(1-t)z)dt\Big|\\
&=&\Big|\sum_{k=1}^{n}(z_{k}-w_{k})\int_{0}^{1}\frac{df}{d\varsigma_{k}}(wt+(1-t)z)dt\\
&&+\sum_{k=1}^{n}(\overline{z}_{k}-\overline{w}_{k})\int_{0}^{1}\frac{df}{d\overline{\varsigma}_{k}}(wt+(1-t)z)dt\Big|\\
&\leq&\sum_{k=1}^{n}|z_{k}-w_{k}|\cdot\Big|\int_{0}^{1}\frac{df}{d\varsigma_{k}}(wt+(1-t)z)dt\Big|\\
&&+\sum_{k=1}^{n}|\overline{z}_{k}-\overline{w}_{k}|\cdot\Big|\int_{0}^{1}\frac{df}{d\overline{\varsigma}_{k}}(wt+(1-t)z)dt\Big|\\
\end{eqnarray*}
\begin{eqnarray*}
&\leq&\Big(\sum_{k=1}^{n}|z_{k}-w_{k}|^{2}\Big)^{\frac{1}{2}}\Big\{\Big[\sum_{k=1}^{n}\Big(\int_{0}^{1}\Big|\frac{\partial
f}{\partial \varsigma_{k}}(wt+(1-t)z)\Big|dt\Big)^{2}\Big]^{\frac{1}{2}}\\
&&+\Big[\sum_{k=1}^{n}\Big(\int_{0}^{1}\Big|\frac{\partial
f}{\partial
\overline{\varsigma}_{k}}(wt+(1-t)z)\Big|dt\Big)^{2}\Big]^{\frac{1}{2}}\Big\}
\\
&\leq& \sqrt{n}|z-w|\Big[\int_{0}^{1}|f_{\varsigma}(wt+(1-t)z)|dt\\
&&+\int_{0}^{1}|f_{\overline{\varsigma}}(wt+(1-t)z)|dt\Big]\\
&\leq& \sqrt{2n}|z-w|\int_{0}^{1}|\widetilde{\nabla f}(wt+(1-t)z)|dt\\
&\leq&C\sqrt{2n}|w-z|\int_{0}^{1}\omega\left(\frac{1}{d(z)-t|w-z|}\right)dt\\
&=&C\sqrt{2n}\int_{0}^{|w-z|}\omega\left(\frac{1}{d(z)-t}\right)dt.
\end{eqnarray*}
This implies
\begin{eqnarray*}
&&\frac{1}{|\mathbb{B}^{n}(z,r)|}\int_{\mathbb{B}^{n}(z,r)}|f(\zeta)-f(z)|dV(\zeta)\\&\leq&
\frac{C\sqrt{2n}}{|\mathbb{B}^{n}(0,r)|}\int_{\mathbb{B}^{n}(0,r)}\left\{\int_{0}^{|\xi|}\omega\left(\frac{1}{d(z)-t}\right)dt\right\}dV(\xi)\\
&=&\frac{C2n\sqrt{2n}}{r^{2n}}\int_{0}^{r}\rho^{2n-1}\left\{\int_{0}^{\rho}\omega\Big(\frac{1}{d(z)-t}\Big)dt\right\}d\rho\\
&\leq&\frac{C2n\sqrt{2n}}{r^{2n}}\int_{0}^{r}\left\{\int_{t}^{r}\rho^{2n-1}
d\rho\right\}\omega\left(\frac{1}{r-t}\right)dt\\
&\leq&\frac{C\sqrt{2n}}{r^{2n}}\int_{0}^{r}(r-t)\left(r^{2n-1}+r^{2n-2}t+\cdots+t^{2n-1}\right)\omega\left(\frac{1}{r-t}\right)dt\\
&\leq&\frac{C\sqrt{2n}}{r^{2n}}r\omega\Big(\frac{1}{r}\Big)\int_{0}^{r}\left(r^{2n-1}+r^{2n-2}t+\cdots+t^{2n-1}\right)dt\\
&=&C\sqrt{2n}\left(\sum_{j=1}^{2n}\frac{1}{j}\right)r\omega\Big(\frac{1}{r}\Big),\\
\end{eqnarray*} where
$\varsigma=(\varsigma_{1},\cdots,\varsigma_{n})=wt+(1-t)z$ and $C$
is a positive constant. The proof of this theorem is complete. \qed

In order to prove Theorem \ref{thm1y}, we need the following lemma.
Using the similar arguments as in the proof of \cite[Lemma 2.5]{MV},
we have

\begin{lem}\label{lem-g1}
Suppose that $f:\ \overline{\mathbb{B}}^{n}(a,r)\rightarrow\mathbb{C}$ is
a continuous function in $\overline{\mathbb{B}}^{n}(a,r)$ and
harmonic  in $\mathbb{B}^{n}(a,r)$. Then
$$|\widetilde{\nabla f}(a)|\leq
\frac{4n\sqrt{n}}{r}\int_{\partial\mathbb{B}^{n}}|f(a+r\zeta)-f(a)|d\sigma(\zeta).$$
\end{lem}
\bpf Let $f=u+iv$, where $u$ and $v$ are real harmonic functions in $\mathbb{B}^{n}(a,r)$. Without loss of generality, we may assume that $a=0$ and $f(0)=0.$
Let
$$K(z,\zeta)=\frac{r^{2n-2}(r^{2}-|z|^{2})}{|z-r\zeta|^{2n}}.$$ Then
$$u(z)=\int_{\partial\mathbb{B}^{n}}K(z,\zeta)u(r\zeta)d\sigma(\zeta),\
z\in\mathbb{B}^{n}(0,r).
$$
By direct calculations, we have
$$\frac{\partial}{\partial x_{j}}K(z,\zeta)=r^{2n-2}
\Big[\frac{-2x_{j}}{|z-r\zeta|^{2n}}-\frac{2n(r^{2}-|z|^{2})(x_{j}-r\alpha_{j})}{|z-r\zeta|^{2n+2}}\Big]$$
and
$$\frac{\partial}{\partial y_{j}}K(z,\zeta)=r^{2n-2}
\Big[\frac{-2y_{j}}{|z-r\zeta|^{2n}}-\frac{2n(r^{2}-|z|^{2})(y_{j}-r\beta_{j})}{|z-r\zeta|^{2n+2}}\Big],$$
which gives
 \be\label{ex}\frac{\partial}{\partial
x_{j}}K(0,\zeta)=\frac{2n\alpha_{j}}{r}~\mbox{and}~\frac{\partial}{\partial
y_{j}}K(0,\zeta)=\frac{2n\beta_{j}}{r}, \ee where
$z=(z_{1},\cdots,z_{n})=(x_{1}+iy_{1},\cdots,x_{n}+iy_{n})$ and
$\zeta=(\alpha_{1}+i\beta_{1},\cdots,\alpha_{n}+i\beta_{n})\in\partial\mathbb{B}^{n}$.
Then by 
(\ref{ex}), we have
\begin{eqnarray*}
|\nabla
u(0)|&=&\Big[\sum_{j=1}^{n}\Big(\Big|\int_{\partial\mathbb{B}^{n}}\frac{\partial}{\partial
x_{j}}K(0,\zeta)u(r\zeta)d\sigma(\zeta)\Big|^{2}\\
&&+\Big|\int_{\partial\mathbb{B}^{n}}\frac{\partial}{\partial
y_{j}}K(0,\zeta)u(r\zeta)d\sigma(\zeta)\Big|^{2}\Big)\Big]^{\frac{1}{2}}\\
&\leq&\sum_{j=1}^{n}\Big(\Big|\int_{\partial\mathbb{B}^{n}}\frac{\partial}{\partial
x_{j}}K(0,\zeta)u(r\zeta)d\sigma(\zeta)\Big|\\
&&+\Big|\int_{\partial\mathbb{B}^{n}}\frac{\partial}{\partial
y_{j}}K(0,\zeta)u(r\zeta)d\sigma(\zeta)\Big|\Big)\\
&\leq&\int_{\partial\mathbb{B}^{n}}|u(r\zeta)|\sum_{j=1}^{n}\Big(\Big|\frac{\partial}{\partial
x_{j}}K(0,\zeta)\Big|+\Big|\frac{\partial}{\partial
y_{j}}K(0,\zeta)\Big|\Big)d\sigma(\zeta)\\
&\leq&\sqrt{2n}\int_{\partial\mathbb{B}^{n}}|u(r\zeta)|\left[\sum_{j=1}^{n}\Big(\Big|\frac{\partial}{\partial
x_{j}}K(0,t)\Big|^{2}+\Big|\frac{\partial}{\partial
y_{j}}K(0,t)\Big|^{2}\Big)\right]^{\frac{1}{2}}d\sigma(\zeta)\\
&=&\frac{2n\sqrt{2n}}{r}\int_{\partial\mathbb{B}^{n}}|u(r\zeta)|d\sigma(\zeta).
\end{eqnarray*}
Similarly, we have
$$|\nabla
v(0)|\leq\frac{2n\sqrt{2n}}{r}\int_{\partial\mathbb{B}^{n}}|v(r\zeta)|d\sigma(\zeta).$$
Then by (\ref{1x}), we conclude that
\begin{eqnarray*} |\widetilde{\nabla f}(0)|&\leq&|f_{z}(0)|+|f_{\overline{z}}(0)|\\
&\leq& |\nabla u(0)|+|\nabla v(0)|\\
&\leq&\frac{2n\sqrt{2n}}{r}\int_{\partial\mathbb{B}^{n}}|u(r\zeta)|+|v(r\zeta)|d\sigma(\zeta)\\
&\leq&\frac{4n\sqrt{n}}{r}\int_{\partial\mathbb{B}^{n}}|f(r\zeta)|d\sigma(\zeta).
\end{eqnarray*}
The proof of the lemma is complete. \epf

\subsection*{Proof of Theorem \ref{thm1y}} First, we show the ``if" part. By Lemma \ref{lem-g1}, we have
$$|\widetilde{\nabla f}(z)|\leq
\frac{4n\sqrt{n}}{\rho}\int_{\partial\mathbb{B}^{n}}|f(z+\rho\zeta)-f(z)|d\sigma(\zeta),$$
where $\rho\in(0,d(z)]$. Let $r=d(z)$. Then we have
$$\int_{0}^{r}|\widetilde{\nabla f}(z)|\rho^{2n} d\rho\leq2\sqrt{n}\int_{0}^{r}\Big(2n\rho^{2n-1}\int_{\partial\mathbb{B}^{n}}|f(z)-f(z+\rho
\zeta)|d\sigma(\zeta)\Big)d\rho,$$ which implies
\begin{eqnarray*}
|\widetilde{\nabla
f}(z)|&\leq&\frac{2(2n+1)\sqrt{n}}{r^{2n+1}}\int_{0}^{r}\Big(2n\rho^{2n-1}\int_{\partial\mathbb{B}^{n}}|f(z)-f(z+\rho
\zeta)|d\sigma(\zeta)\Big)d\rho\\
&=&\frac{2(2n+1)\sqrt{n}}{r|\mathbb{B}^{n}(z,r)|}\int_{\mathbb{B}^{n}(z,r)}|f(\xi)-f(z)|dV(\xi)\\
&\leq&2(2n+1)\sqrt{n}C\omega\Big(\frac{1}{r}\Big)\\
&=&2(2n+1)\sqrt{n}C\omega\left(\frac{1}{d(z)}\right).
\end{eqnarray*}

The ``only if" part easily follows from Theorem \ref{thm1-y}. The
proof of the theorem is complete. \qed

\section{The finite Dirichlet energy integral  and its application}\label{csw-sec3}

\begin{lem}\label{thm2}
Let $f$ be a solution to \eqref{eq-g2}. Then for $p\in[2,\infty)$
and $\beta\in(0,\infty)$,
$$D_{f}(\beta,p-2,2)\leq \frac{\beta\sqrt{2}}{2} D_{f}(\beta-1,p-1,1).$$
\end{lem}

\bpf By  Lemmas \ref{lem1} and \ref{Lemma-2.0}, we have
\beq\label{eq-1-thm2} \nonumber
r^{2n-1}\frac{d}{dr}M_{p}^{p}(r,f)&=&\frac{1}{2n}\int_{\mathbb{B}^{n}(0,r)}
\Big[p(p-2)|f(z)|^{p-4}\sum_{k=1}^{n}|f(z)\overline{f_{z_{k}}(z)}+\overline{f(z)}f_{\overline{z}_{k}}(z)|^{2}\\
 &&+2p|f(z)|^{p-2}|\widetilde{\nabla
f}(z)|^{2}+p\lambda|f(z)|^{p}\Big]dV_{N}(z), \eeq which implies
\beq\label{eq-2-thm2}
  \nonumber \frac{d}{dr}\left(r^{2n-1}\frac{d}{dr}M_{p}^{p}(r,f)\right)&=&\int_{\partial\mathbb{B}^{n}}
  r^{2n-1}\Big[2p|f(r\zeta)|^{p-2}|\widetilde{\nabla
f}(r\zeta)|^{2}\\   \nonumber &&+
p(p-2)|f(r\zeta)|^{p-4}\sum_{k=1}^{n}|f(r\zeta)\overline{f_{z_{k}}(r\zeta)}
+\overline{f(r\zeta)}f_{\overline{z}_{k}}(r\zeta)|^{2}\\
&&+p\lambda|f(r\zeta)|^{p}\Big] d\sigma(\zeta).\eeq

In addition, we see \beq\label{eqt-1}\nonumber
\frac{d}{dr}M_{p}^{p}(r,f)&=&
\int_{\partial\mathbb{B}^{n}}\frac{d}{dr}(|f(r\zeta)|^{p})d\sigma(\zeta)\\
\nonumber &=&p\int_{\partial\mathbb{B}^{n}}|f(r\zeta)|^{p-2}
\mbox{Re}\left[\sum_{k=1}^{n}\Big(f_{z_{k}}(r\zeta)\overline{f(r\zeta)}+f(r\zeta)
\overline{f_{\overline{z}_{k}}(r\zeta)}\Big)\zeta_{k}\right]d\sigma(\zeta)\\
 &\leq&
p\sqrt{2}\int_{\partial\mathbb{B}^{n}}|f(r\zeta)|^{p-1}|\widetilde{\nabla
f}(r\zeta)|d\sigma(\zeta),\eeq where
$\zeta=(\zeta_{1},\cdots,\zeta_{n})\in\partial\mathbb{B}^{n}$.

It follows from (\ref{eq-2-thm2}) and (\ref{eqt-1}) that
 \beq\label{eq-3-thm2}
&&\beta
p\sqrt{2}\int_{\mathbb{B}^{n}}(1-|z|)^{\beta-1}|f(z)|^{p-1}|\widetilde{\nabla
f}(z)|dV_{N}(z)\\ \nonumber
 &=& \beta p\sqrt{2}\int_{0}^{1}2nr^{2n-1}(1-r)^{\beta-1}\int_{\partial\mathbb{B}^{n}}|f(r\zeta)|^{p-1}|\widetilde{\nabla
f}(r\zeta)|d\sigma(\zeta) dr\\ \nonumber&\geq&
\beta\int_{0}^{1}2nr^{2n-1}(1-r)^{\beta-1}\left(\frac{d}{dr}M_{p}^{p}(r,f)\right)dr\\
\nonumber
 &=&
 \int_{0}^{1}2n(1-r)^{\beta}\frac{d}{dr}\left(r^{2n-1}\frac{d}{dr}M_{p}^{p}(r,f)\right)dr\\ \nonumber
 \eeq
\beq
 \nonumber
 &=&\int_{0}^{1}2n(1-r)^{\beta}\int_{\partial\mathbb{B}^{n}}
  r^{2n-1}\Big[2p|f(r\zeta)|^{p-2}|\widetilde{\nabla
f}(r\zeta)|^{2}\\   \nonumber &&+
p(p-2)|f(r\zeta)|^{p-4}\sum_{k=1}^{n}|f(r\zeta)\overline{f_{z_{k}}(r\zeta)}
+\overline{f(r\zeta)}f_{\overline{z}_{k}}(r\zeta)|^{2}\\ \nonumber
&&+p\lambda|f(r\zeta)|^{p}\Big] d\sigma(\zeta)dr\\ \nonumber
 &=& 2p\int_{\mathbb{B}^{n}}\Big[|f(z)|^{p-2}|\widetilde{\nabla
 f}(z)|^{2}+\frac{\lambda}{2}|f(z)|^{p}+\\ \nonumber
 &&\big(\frac{p}{2}-1\big)|f(z)|^{p-4}
\sum_{k=1}^{n}\left|f(z)\overline{f_{z_{k}}(z)}+\overline{f(z)}f_{\overline{z}_{k}}(z)\right|^{2}\Big](1-|z|)^{\beta}dV_{N}(z)\\
\nonumber &
\geq&2p\int_{\mathbb{B}^{n}}|f(z)|^{p-2}|\widetilde{\nabla
 f}(z)|^{2}(1-|z|)^{\beta}dV_{N}(z), \eeq
whence
$$D_{f}(\beta,p-2,2)\leq \frac{\beta\sqrt{2}}{2} D_{f}(\beta-1,p-1,1),$$
from which the proof follows. \epf

By  elementary computations, we easily see that
\begin{lem}\label{Lem1}
Suppose that $a,~b\in[0,\infty)$ and $q\in(0,\infty)$. Then
$$(a+b)^{q}\leq2^{\max\{q-1,0\}}(a^{q}+b^{q}).$$
\end{lem}

\subsection*{Proof of  Theorem \ref{thm4}} By Lemma \ref{lem1}, we know that $|\widetilde{\nabla f}|^{2}$ is
subharmonic in $\mathbb{B}^{n}$. Then for $r\in[0,1-|z|)$, we have
$$|\widetilde{\nabla f}(z)|^{2}\leq\int_{\partial\mathbb{B}^{n}}|\widetilde{\nabla
f}(z+r\zeta)|^{2}d\sigma(\zeta).$$ Integration and Lemma \ref{thm2}
yield
 \begin{eqnarray*}
\frac{(1-|z|)^{2n}|\widetilde{\nabla f}(z)|^{2}}{2^{2n}}
&\leq&\int_{\partial\mathbb{B}^{n}}\int_{0}^{\frac{1-|z|}{2}}2nr^{2n-1}|\widetilde{\nabla
f}(z+r\zeta)|^{2}drd\sigma(\zeta)\\
&=&\int_{\mathbb{B}^{n}(z,\frac{1-|z|}{2})}|\widetilde{\nabla f}(\xi)|^{2}dV_{N}(\xi)\\
&\leq&2^{\beta}(1-|z|)^{-\beta}\int_{\mathbb{B}^{n}(z,\frac{1-|z|}{2})}(1-|\xi|)^{\beta}|\widetilde{\nabla
f}(\xi)|^{2}dV_{N}(\xi)\\
&\leq&2^{\beta}D_{f}(\beta,0,2)(1-|z|)^{-\beta}\\
 &\leq&\beta2^{\beta-\frac{1}{2}}
D_{f}(\beta-1,1,1)(1-|z|)^{-\beta}
\end{eqnarray*}
which gives \be\label{eq-1}|\widetilde{\nabla
f}(z)|\leq\frac{C_{3}}{(1-|z|)^{n+\frac{\beta}{2}}},\ee where
$C_{3}=\sqrt{\beta2^{\beta-1/2+2n}D_{f}(\beta-1,1,1)}$.

By (\ref{eq-1}), we have
\begin{eqnarray*}
|f(z)|&\leq& |f(0)|+\left|\int_{[0,z]}df(\zeta)\right|\\
&\leq&|f(0)|+\sqrt{2}\int_{[0,z]}|\widetilde{\nabla f}(\zeta)||d\zeta|\\
&\leq&|f(0)|+\frac{C_{4}}{(1-|z|)^{\frac{\beta}{2}+n-1}},
\end{eqnarray*}
where $C_{4}=\sqrt{2}C_{3}/(n-1+\beta/2)$ and $[0,z]$ denotes the
segment from $0$ to $z$. Then by Lemma \ref{Lem1}, we see that for
$z\in\mathbb{B}^{n}$,
\beq\label{eq-2}|f(z)|^{\frac{2}{\beta}}&\leq&\left[|f(0)|+\frac{C_{4}}{(1-|z|)^{\frac{\beta}{2}+n-1}}\right]^{\frac{2}{\beta}}\\
\nonumber
&\leq&2^{\frac{2}{\beta}-1}\left[|f(0)|^{\frac{2}{\beta}}+\frac{C_{4}^{\frac{2}{\beta}}}{(1-|z|)^{1+\frac{2(n-1)}{\beta}}}\right]\eeq
and
\beq\label{eq-3}|f(z)|^{\frac{2}{\beta}-2}&\leq&\left[|f(0)|+\frac{C_{4}}{(1-|z|)^{\beta/2+n-1}}\right]^{\frac{2}{\beta}-2}\\
\nonumber
&\leq&2^{\frac{2}{\beta}-2}\left[|f(0)|^{\frac{2}{\beta}-2}+
\frac{C_{4}^{\frac{2}{\beta}-2}}{(1-|z|)^{1-\beta+(n-1)(\frac{2}{\beta}-2)}}\right].\eeq

Let $p=2/\beta.$ We divide the rest of the proof into two cases.
\bca Let $p\in[4,\infty)$.\eca By direct calculations, we get
 \beq\label{eq-4} \nonumber
\Delta(|f|^{p})&=&4\sum_{k=1}^{n}\frac{\partial^{2}}{\partial
z_{k}\partial\overline{z}_{k}}(|f|^{p})\\
\nonumber&=&p(p-2)|f|^{p-4}\sum_{k=1}^{n}|f_{z_{k}}\overline{f}+\overline{f}_{\overline{z}_{k}}f|^{2}+2p|f|^{p-2}|\widetilde{\nabla
f}|^{2}+p\lambda|f|^{p}\\ &\leq&2p(p-1)|f|^{p-2}|\widetilde{\nabla
f}|^{2}+p\lambda|f|^{p}.
 \eeq
Hence by  (\ref{eq-2}), (\ref{eq-3}) and (\ref{eq-4}), we conclude
that for $z\in\mathbb{B}^{n}$, \beq\label{eq-5}\nonumber
(1-|z|)^{1+p(n-1)}\Delta(|f(z)|^{p})&\leq&2p(p-1)(1-|z|)^{1+p(n-1)}|f(z)|^{p-2}|\widetilde{\nabla
f}(z)|^{2}\\ \nonumber &&+p\lambda(1-|z|)^{1+p(n-1)}|f(z)|^{p}\\
\nonumber &\leq&2p(p-1)(1-|z|)^{\beta}|\widetilde{\nabla
f}(z)|^{2}(1-|z|)^{1+p(n-1)-\beta}|f(z)|^{p-2}\\ \nonumber
&&+p\lambda2^{p-1}\big(C_{4}^{p}+|f(0)|^{p}\big)\\
&\leq&C_{5}+C_{6}(1-|z|)^{\beta}|\widetilde{\nabla f}(z)|^{2}, \eeq
 where
$C_{5}=p\lambda2^{p-1}\big(C_{4}^{p}+|f(0)|^{p}\big)$ and
$C_{6}=2p(p-1)2^{p-2}\big(|f(0)|^{p-2}+C_{4}^{p-2}\big)$. By Theorem
\ref{thm2}, we know
\be\label{eq1-th-1}D_{f}(\beta,0,2)\leq\frac{\beta\sqrt{2}}{2}
D_{f}(\beta-1,1,1).\ee Therefore, (\ref{eq-5}) and (\ref{eq1-th-1})
imply that there exist positive constants $C_{1}$ and $C_{2}$ such
that
$$\int_{\mathbb{B}^{n}}(1-|z|)^{1+p(n-1)}\Delta(|f(z)|^{p})dV_{N}(z)\leq C_{1}D_{f}(\frac{2}{p}-1,1,1)+C_{2}.$$

\bca Let $p\in[2,4)$.\eca In this case, we let
$F_{m}^{p}=(|f|^{2}+\frac{1}{m})^{p/2}$, and let
$T_{m}=\Delta(F_{m}^{p})$. Obviously, for   $r\in(0,1)$, $T_{m}$ is
integrable in $\mathbb{B}^{n}(0,r)$ and $T_{m}\leq F,$ where
$$F=p(p-2)|f|^{p-2}\sum_{k=1}^{n}(|f_{z_{k}}|+|f_{\overline{z}_{k}}|)^{2}+2p(1+|f|^{2})^{\frac{p}{2}-1}|\widetilde{\nabla
f}|^{2}+p\lambda|f|^{2}(|f|^{2}+1)^{\frac{p}{2}-1}$$ and $F$ is
integrable in $\mathbb{B}^{n}(0,r)$.

Then, by Lebesgue's Dominated Convergence Theorem together with (\ref{eq-5}), 
we have
\beq\label{eq-t}\nonumber
&&\lim_{n\rightarrow\infty}\int_{\mathbb{B}^{n}(0,r)}(1-|z|)^{1+p(n-1)}\Delta(F_{m}^{p}(z))dV_{N}(z)\\
\nonumber &=&\int_{\mathbb{B}^{n}(0,r)}(1-|z|)^{1+p(n-1)}
\lim_{n\rightarrow\infty}\big[\Delta(F_{m}^{p}(z))\big]dV_{N}(z)\\
\nonumber &=&p\int_{\mathbb{B}^{n}(0,r)}\Big [\big (p-2\big
)|f(z)|^{p-4}
\sum_{k=1}^{n}|f_{z_{k}}(z)\overline{f(z)}+\overline{f_{\overline{z}_{k}}(z)}f(z)|^{2}\\
\nonumber
 &&  +2|f(z)|^{p-2}|\widetilde{\nabla
f}(z)|^{2}+\lambda|f(z)|^{p}\Big](1-|z|)^{1+p(n-1)}\, dV_{N}(z)\\
\nonumber
&\leq&\int_{\mathbb{B}^{n}(0,r)}[C_{5}+C_{6}(1-|z|)^{\beta}|\widetilde{\nabla
f}(z)|^{2}]dV_{N}(z),\eeq and so we infer from    (\ref{eq-2}),
(\ref{eq-3}) and Theorem \ref{thm2} that there exist positive
constants $C_{1}$ and $C_{2}$ such that
$$\int_{\mathbb{B}^{n}}(1-|z|)^{1+p(n-1)}\Delta(|f(z)|^{p})dV_{N}(z)\leq C_{1}D_{f}(\frac{2}{p}-1,1,1)+C_{2}.$$
 The proof of this theorem is complete. \qed

\subsection*{Proof of  Corollary \ref{cor-1}}
 For a fixed $r\in(0,1)$,
since $$\lim_{|z|\rightarrow r}\frac{\log r-\log
|z|}{r-|z|}=\frac{1}{r},$$  we see that there is $r_{0}\in(0,r)$
satisfying
$$\log r-\log |z|\leq\frac{2}{r}(r-|z|)$$ for $r_{0}\leq|z|<r$. Let
$p=2/\beta$. The it follows from
$$\lim_{\rho\rightarrow0+}\rho\log\frac{1}{\rho}=0$$  that
 \beq\label{eq-x5} \int_{\mathbb{D}_{r_{0}}}\Delta
(|f(z)|^{p})\log\frac{r}{|z|}\,d\sigma(z)&\leq&\int_{\mathbb{D}_{r_{0}}}\Delta
(|f(z)|^{p})\log\frac{1}{|z|}\,d\sigma(z)\\ \nonumber
&=&\int_{0}^{2\pi}\int_{0}^{r_{0}}\Delta (|f(\rho
e^{i\theta})|^{p})\rho\log\frac{1}{\rho}d\rho d\theta\\
\nonumber&<&\infty. \eeq Since $D_{f}(\beta-1,1,1)<\infty$, it
follows from Theorem \ref{thm4} that
\be\label{cor-2v}\int_{\mathbb{D}\backslash
\mathbb{D}_{r_{0}}}\Delta (|f(z)|^{p})(1-|z|)\,d\sigma(z)<\infty.\ee

 Hence by (\ref{eq-x5}), (\ref{cor-2v}) and Theorem
\Ref{Green-thm}, we obtain that
\begin{eqnarray*}
M_{p}^{p}(r,f)&=&|f(0)|^{p}+ \frac{1}{2}\int_{\mathbb{D}_{r}}\Delta
(|f(z)|^{p})\log\frac{r}{|z|}\,d\sigma(z)\\
&=&|f(0)|^{p}+\frac{1}{2}\int_{\mathbb{D}_{r_{0}}}\Delta
(|f(z)|^{p})\log\frac{r}{|z|}\,d\sigma(z)\\
&&+\frac{1}{2}\int_{\mathbb{D}_{r}\backslash
\mathbb{D}_{r_{0}}}\Delta
(|f(z)|^{p})\log\frac{r}{|z|}\,d\sigma(z)\\
 &\leq&|f(0)|^{p}+\frac{1}{2}\int_{\mathbb{D}_{r_{0}}}\Delta
(|f(z)|^{p})\log\frac{r}{|z|}\,d\sigma(z)\\
&&+ \int_{\mathbb{D}_{r}\backslash \mathbb{D}_{r_{0}}}\Delta
(|f(z)|^{p})\frac{(r-|z|)}{r}\,d\sigma(z)\\
&\leq&|f(0)|^{p}+\frac{1}{2}\int_{\mathbb{D}_{r_{0}}}\Delta
(|f(z)|^{p})\log\frac{r}{|z|}\,d\sigma(z)\\
&&+ \int_{\mathbb{D}\backslash \mathbb{D}_{r_{0}}}\Delta
(|f(z)|^{p})(1-|z|)\,d\sigma(z)\\
 &<&\infty.
 \end{eqnarray*}
Since Lemma \ref{lem1} shows that the function $M_{p}^{p}(r,f)$ is
increasing with respect to $r$ in $(0,1)$, we know that the limit
$$\lim_{r\rightarrow1-}M_{p}(r,f)$$ does exist, which implies
$f\in\mathcal{H}^{p}$. The proof of the corollary is complete.
 \qed


\end{document}